\documentclass[11pt]{article}

\usepackage{amssymb}
\usepackage{amsmath}
\usepackage{amsthm}
\usepackage{graphicx}
\usepackage{siunitx}
\usepackage{subfig}

\newcommand{\bg}{{\boldsymbol g}}
\newcommand{\bG}{{\boldsymbol G}}
\newcommand{\be}{{\boldsymbol e}}

\newcommand{\bm}{{\boldsymbol m}}
\newcommand{\bn}{{\boldsymbol n}}
\newcommand{\bu}{{\boldsymbol u}}
\newcommand{\bhu}{{\boldsymbol{\hat{u}}}}
\newcommand{\bv}{{\boldsymbol v}}
\newcommand{\bw}{{\boldsymbol w}}
\newcommand{\bx}{{\boldsymbol x}}
\newcommand{\bbeta}{{\boldsymbol \beta}}
\newcommand{\bxi}{{\boldsymbol \xi}}
\newcommand{\bpsi}{{\boldsymbol \psi}}
\newcommand{\bpi}{{\boldsymbol \pi}}
\newcommand{\bchi}{{\boldsymbol \chi}}
\newcommand{\bphi}{{\boldsymbol \phi}}
\newcommand{\wdot}{\dot{w}}
\newcommand{\btau}{{\boldsymbol \tau}}
\newcommand{\oM}{\overline{M}}
\newcommand{\bitv}{\boldsymbol{\mathit{v}}}
\newcommand{\ds}{\displaystyle}
\newcommand{\R}{\mathbb R}
\newcommand{\Po}{\mathbb P}
\newcommand{\cT}{{\cal T}}

\author{Z.~Sun, M.~Braack, J.~Lang}
\title{An Adaptive Moving Finite Element Method for
Steady Low Mach Number Compressible Combustion Problems}
\author{
Zhen Sun \\
{\small \it Technische Universit\"at Darmstadt} \\
{\small \it Dolivostra{\ss}e 15, 64293 Darmstadt, Germany} \\
{\small sun@mathematik.tu-darmstadt.de} \\ \\
Malte Braack \\
{\small \it Christian-Albrechts-Universit\"at zu Kiel} \\
{\small \it Westring 383, 24098 Kiel, Germany} \\
{\small braack@math.uni-kiel.de} \\ \\
Jens Lang\footnote{corresponding author}\\
{\small \it Technische Universit\"at Darmstadt} \\ {\small \it
Dolivostra{\ss}e 15, 64293 Darmstadt, Germany}\\
{\small lang@mathematik.tu-darmstadt.de}}
\date{March 26, 2019}

\begin{document}
\maketitle

\begin{abstract}
This work surveys an $r$-adaptive moving mesh finite element method for the numerical solution of premixed laminar flame problems. Since the model of chemically reacting
flow involves many different modes with diverse length scales, the computation of such a problem is often extremely time-consuming. Importantly, to capture the significant characteristics of the flame structure when using detailed chemistry, a much more stringent requirement on the spatial resolution of the interior layers of some intermediate species is necessary. Here, we propose a moving mesh method in which the mesh is obtained from the solution of so-called moving mesh partial differential equations. Such equations result from the variational formulation of a minimization problem for a given target functional that characterizes the inherent difficulty in the numerical approximation of the underlying physical equations. Adaptive mesh movement has emerged as an area of intense research in mesh adaptation in the last decade. With this approach points are only allowed to be shifted in space leaving the topology of the grid unchanged. In contrast to methods with local refinement, data structure hence is unchanged and load balancing is not an issue as grid points remain on the processor where they are. We will demonstrate the high potential of moving mesh methods for effectively optimizing the distribution of grid points to reach the required resolution for chemically reacting flows with extremely thin boundary layers.
\end{abstract}


\section{Introduction}
The numerical simulation of chemically reacting flows with detailed
chemistry is still a challenging task due to the presence of
a large range of multi-scale aspects. High spatial resolution is
needed in the vicinity of the flame zone, while simultaneously
coarser meshes can be used in regions with relatively large flow
structures, e.g., downstream of the flame. Adaptive mesh refinement (AMR)
based on dynamically refining or coarsening the computational mesh
to match local features has proven to be a very efficient strategy to
perform a multi-scale combustion simulation. An overview of the basic
design concepts used to develop block-structured AMR algorithms
for the low Mach number model has been given by {\sc Bell} and {\sc Day} \cite{BellDay2011}. Two of the authors, {\sc Braack} and {\sc Lang}, have also been contributing to the development of AMR methods \cite{BeckerBraackRannacher1999,BraackRichter2006,FroehlichLang1998,Lang1996}.
However, since reacting flow applications with detailed
chemistry and transport models already tend to be significantly complex, an AMR implementation on parallel computers might quickly become prohibitively intricate.
In such a situation, adaptive moving mesh methods are considered to be an attractive
alternative since the mesh topology and hence the underlying data structures
are left unchanged. Such methods have become available through the implementation of variational principles that define a moving mesh as the solution of so-called moving mesh partial differential equations (MMPDE) introduced by {\sc Huang} and
{\sc Russell}, see \cite{HuangRussell2010} for a thorough overview.
Recently, adaptive moving meshes (AMM) have been successfully applied
to large eddy simulations of complex turbulent flows by {\sc Liersch}, {\sc Frankenbach},
{\sc Fr\"ohlich}, and {\sc Lang} \cite{LierschFrankenbachFroehlichLang2014}.

In this paper, we investigate the application of the AMM technique to steady low
Mach number compressible combustion problems in two spatial dimensions. Special emphasize
is put on the design of the monitor function that constitutes the heart of the MMPDE and links
the movement of the grid points to the physical solution. We will discuss two approaches that
are based on the gradient and the curvature of selected chemical components, respectively. The latter one is motivated by a standard error estimate of the interpolation error. The steady combustion model
is discretized by stabilized linear finite elements that use pressure stabilization and
streamline diffusion to enhance the main diagonal of the resulting algebraic system. A
pseudo linearly implicit time-stepping scheme is applied to solve the highly nonlinear algebraic
systems. Two examples are considered to demonstrate the potential of our AMM method: a benchmark three-component ozone decomposition and a methane flame in the practical application of a
prototype lamelle burner from the {\sc Bosch} company.

\section{Theoretical Fundamentals of Reactive Flows}
\label{sec:modeling}
\subsection{Models for stationary chemically reacting flows}
We denote the flow velocities by $\bv$, the pressure by $p$, the temperature by $T$ and the density by $\rho$. Furthermore, $\bw=(w_1,\ldots,w_{n_s})^T$ is defined as the
vector of mass fractions $w_i$ of the $i$-th species, where $i=1,\ldots,n_s$. The basic equations for a stationary reactive viscous flow express the conservation and balance laws for the total mass, momentum, energy and each species, respectively, accomplished by the equation of state:
\begin{eqnarray}
\label{eq:viscous_euler_mass}
\nabla \cdot (\rho \bv) &=& 0, \\[1mm]
\label{eq:viscous_euler_momentum}
\rho ( \bv \cdot \nabla ) \bv + \nabla p + \nabla \cdot \btau &=& \rho \bg, \\[1mm]
\rho \bv \cdot \nabla w_i + \nabla \cdot ( \rho D_i \nabla w_i ) &=& M_i \,\wdot_i(T, \boldsymbol \omega), \quad i=1,\ldots, n_s,\\[1mm]
\label{eq:viscous_euler_temperature}
c_p \rho \bv \cdot \nabla T - \nabla \cdot (\lambda \nabla T)
- \bv \cdot \nabla p + \btau : \nabla \bv &=& -\sum_{i=1}^{n_s} h_i M_i
\wdot_i (T, \bw),\\[1mm]
\label{eq:viscous_euler_eos}
p &=& \frac{\rho RT}{\oM},\quad
\oM = \left( \sum_{i=1}^{n_s} \frac{w_i}{M_i} \right)^{-1},
\end{eqnarray}
where $\bg$ is the gravitational force and $c_p$ is the heat capacity of the mixture at constant pressure. For each species, $M_i$ is the molecular weight, $h_i$ its (mass-)specific enthalpy, $\wdot_i$ its molar production rate, and $D_i$ its mass
diffusion coefficient. The viscous stress tensor $\btau$ is given by
\begin{equation}
\btau = - \mu  \left(\nabla \bv + (\nabla \bv)^T
- \frac{2}{3}(\nabla \cdot \bv) \textbf{I} \right)
\end{equation}
with $\mu$ being the dynamic viscosity of the fluid.
Compressible flow equations admit material and acoustic waves that propagate
at velocity $\bv$ and speed of sound $c$, respectively. In most practical
combustion systems, the flow is in the low Mach number regime, i.e., $M=|\bv|/c\ll 1$.
This disparity in scales can be exploited to compute combustion problems much
more efficiently. A rigorous low Mach number asymptotic analysis \cite{MadjaSethian1985}
for the behaviour $M\rightarrow 0$ shows that the total pressure can be decomposed as
\begin{equation}
 p(\bx)  = P_{th} + p_{hyd} (\bx),
\end{equation}
where $P_{th}$ is the constant thermodynamic pressure and
$p_{hyd}(\bx)$ is the hydrodynamic pressure. With this decomposition, $P_{th}$ defines
the thermodynamic state and $p_{hyd}$ plays the role of a Lagrange multiplier to constrain
the flow so that the thermodynamic pressure is equilibrated everywhere in space. Hence,
compressibility effects due to chemical heat release and other thermal processes are
retained while acoustic wave propagation is entirely eliminated.
Eventually, the simplified low Mach number approximation for steady compressible reactive flows becomes \cite{BraackRannacher1999}
\begin{eqnarray}
\label{eq:low_Mach_number_mass}
\nabla \cdot \bv &=& \frac{1}{T} \bv \cdot \nabla T -
\frac{1}{\oM} \bv \cdot \nabla \oM, \\[1mm]
\label{eq:low_Mach_number_momentum}
\rho ( \bv \cdot \nabla ) \bv + \nabla p_{hyd} -
\nabla \cdot ( \mu \nabla \bv ) &=& \rho \bg, \\[1mm]
\rho \bv \cdot \nabla w_i + \nabla \cdot ( \rho D_i \nabla w_i ) &=& M_i \,
\wdot_i(T,\bw), \quad i=1,\ldots, n_s,\\[1mm]
\label{eq:low_Mach_number_temperature}
c_p \rho \bv \cdot \nabla T - \nabla \cdot (\lambda \nabla T)
&=& -\sum_{i=1}^{n_s} h_i M_i
\wdot_i (T, \bw),\\[1mm]
\label{eq:low_Mach_number_eos}
\rho &=& \frac{P_{th}\oM}{RT},\quad
\oM = \left( \sum_{i=1}^{n_s} \frac{w_i}{M_i} \right)^{-1}.
\end{eqnarray}
Here, we have differentiated the new equation of state in
(\ref{eq:low_Mach_number_eos}) to derive a
constraint on the velocity in (\ref{eq:low_Mach_number_mass}). As
usual, parts of the stress tensor have been absorbed into $p_{hyd}$.
We will only consider $n_s-1$ species
equations together with the requirement $\sum_{i=1}^{n_s}w_i=1$ to
keep mass conservation and the computation of
diffusion velocities consistent \cite{SmookeGiovangigli1991}. Since
only laminar flames are considered, the influence of the gravitation
in (\ref{eq:low_Mach_number_momentum}) is also often neglected. We set
$P_{th}=101\,325\,\si{Pascal}$ in our applications.

\subsection{Finite Element Discretization}
To simplify the notation and write the system
(\ref{eq:low_Mach_number_mass})-(\ref{eq:low_Mach_number_eos}) in a more compact form,
we extend the vector of mass fractions by $w_0:=T$ and introduce
\begin{equation}
\nu_0 = \lambda, \quad \nu_i = \rho D_i, \;i=1,\ldots,n_s,
\end{equation}
for the diffusion coefficients of all components of $\bw=(w_0,\ldots,w_{n_s})^T$.
With the new quantities
\begin{align}
&\;\bbeta = \rho \bv,\quad\bm=\oM^{-1}\nabla \oM,\quad l=(\rho T)^{-1},\\
&\;f_0(\bw)= -\sum_{i=1}^{n_s} h_i M_i\wdot_i,\quad
f_i(\bw)=M_i\,\wdot_i,\;i=1,\ldots,n_s,
\end{align}
and the convention that $\bar{\bbeta}=c_p\bbeta$ for the temperature equation
($i=0$) and $\bar{\bbeta}=\bbeta$ otherwise,
the stationary low Mach number system can be rewritten in the condensed form
\begin{align}
\label{eq:lM_compact0}
\nabla \cdot \bv - l \bbeta \cdot \nabla w_0 + \bv \cdot \bm &\;=\; 0, \\[1mm]
\label{eq:lM_compact1}
( \bbeta \cdot \nabla ) \bv - \nabla \cdot (\mu \nabla \bv ) + \nabla p_{hyd} &\;=\;
\rho \bg, \\[1mm]
\label{eq:lM_compact2}
\bar{\bbeta} \cdot \nabla w_i - \nabla \cdot ( \nu_i \nabla w_i ) &\;=\; f_i(\bw),\;
i=0,\ldots,n_s.
\end{align}
For the presentation of the finite element approximation of this system, we follow
the approach presented in \cite[Sec. 2]{BraackRannacher1999}, see also \cite{Braack1998}.

Let $\Omega\in\R^2$ be our bounded computational domain with polygonal boundary
$\partial\Omega$. Introducing the inner product and norm
\begin{equation}
(f,g) := \int_\Omega f(\bx)g(\bx)\,d\bx, \quad
\Vert f \Vert := \left( \int_\Omega \vert f(\bx) \vert^2 \,d\bx \right)^{1/2},
\end{equation}
the Lebesgue space of all square-integrable functions on $\Omega$, $L^2(\Omega)$,
is defined by all functions $f(\bx)$ with $\Vert f \Vert < \infty$. The functions
from $L^2(\Omega)$ with square-integrable generalized first-order derivatives form
the Sobolev space $H^1(\Omega)$. We will use certain subspaces of these two spaces to
set up a variational weak formulation for (\ref{eq:lM_compact0})-(\ref{eq:lM_compact2}),
namely
\begin{equation}
p\in Q=L^2(\Omega)/\R,\quad \bv\in H\subset \left( H^1(\Omega)\right)^2,
\quad \bw\in R\subset \left(H^1(\Omega)\right)^{n_s+1}.
\end{equation}
Note that the pressure is determined only modulo constants, which is expressed
by the special construction of the space $Q$. The solution triple $\bu:=(p,\bv,\bw)$
is now an element of the space $V:=Q\times H\times R$. Testing the equations
(\ref{eq:lM_compact0})-(\ref{eq:lM_compact2}) with a function
$\bphi=(\theta,\bchi,\bpi)\in V$, integrating over $\Omega$ and applying integration
by parts for the diffusive terms and the pressure gradient yields the variational
nonlinear equations
\begin{equation}
\label{eq:var_form}
\bu\in V: \quad A(\bu,\bphi) = 0 \quad \text{for all } \bphi\in V
\end{equation}
with the semi-linear form
\begin{align}
A(\bu,\bphi) := &\; ( \nabla \cdot \bv, \theta ) - (l \bbeta \cdot \nabla w_0, \theta )
+ ( \bv \cdot \bm, \theta ) \,+ \nonumber \\[1mm]
&\; (( \bbeta \cdot \nabla ) \bv, \bchi ) + ( \mu \nabla \bv, \nabla \bchi )
- ( p_{hyd}, \nabla \cdot \bchi ) - ( \rho \bg, \bchi ) \,+ \nonumber \\[1mm]
&\; \ds \sum_{i=0}^{n_s} \lbrace ( \bar{\bbeta} \cdot \nabla w_i, \pi_i ) +
(\nu_i \nabla w_i, \nabla \pi_i ) - ( f_i(\bw),\pi_i) \rbrace.
\end{align}
Here, for simplicity, we have used the free-stream outflow condition and homogeneous
Dirichlet conditions. Other boundary conditions can be handled by natural modifications.

We will now describe our finite element approximation. First, we decompose $\overline{\Omega}$
into a regular partition of triangles $\cT_h=\{K\}$ and
define the conforming finite element space of continuous piecewise linear functions by
\begin{equation}
S_h := \{ f \in C^0(\overline{\Omega}):\;f_{|K}\in\Po_1,\;K\in \cT_h \},
\end{equation}
where $\Po_1$ is the space of all polynomials of degree not larger than one.
Then, the infinite dimensional space $V$ is approximated by a finite dimensional space
$V_h=(Q_h,H_h,R_h)$ and the corresponding discrete approximations $\bu_h=(p_h,\bv_h,\bw_h)$ are
determined in the finite element spaces
\begin{equation}
Q_h := S_h/\R, \quad H_h := \left(S_h\right)^2,\quad R_h := \left(S_h\right)^{n_s+1}.
\end{equation}
The main observation is now that a simple replacement of $V$ by $V_h$ in (\ref{eq:var_form})
does not give a stable discretization. Instead, we use pressure stabilization and
streamline diffusion to enhance the main diagonal of the resulting algebraic system.
This gives the finite element approximation
\begin{equation}
\label{eq:var_form_fem}
\bu_h\in V_h: \quad A(\bu_h,\bphi) + S_h(\bu_h,\bphi)= 0 \quad \text{for all } \bphi\in V_h
\end{equation}
with the stabilization term $S_h$ defined by
\begin{equation}
S_h(\bu_h,\bphi) := c_h(\bu_h,\theta) + m_h(\bu_h,\bchi) + \sum_{i=0}^{n_s} t_{i,h}(\bu_h,\pi_i),
\end{equation}
where
\begin{align}
c_h(\bu_h,\theta) &\;:=\; \sum_{K\in\cT_h} \alpha_K \left(
( \bbeta \cdot \nabla ) \bv + \nabla p_{hyd} - \rho \bg, \nabla \theta )\right)_K,\\[2mm]
m_h(\bu_h,\bchi)  &\;:=\; \sum_{K\in\cT_h} \delta_K \left(
( \bbeta \cdot \nabla ) \bv + \nabla p_{hyd} - \rho \bg, \bbeta\cdot\nabla\bchi )\right)_K,\\[2mm]
t_{i,h}(\bu_h,\pi_i) &\;:=\; \sum_{K\in\cT_h} \delta_{K,i} \left(
( \bar{\bbeta} \cdot \nabla ) w_{h,i} - f_i(\bw_h), \bar{\bbeta}\cdot\nabla\pi_i )\right)_K,\;i=0,\ldots,n_s.
\end{align}
Here, $(\cdot,\cdot)_K$ stands for the inner product over the triangle $K$.
The density $\rho_h$ is a function of $\bw_h$ and determined from the equation of state in
(\ref{eq:low_Mach_number_eos}). The mesh-dependent constants $\alpha_k$, $\beta_K$, and
$\beta_{K,i}$ must be chosen carefully in order to locally reflect convection-dominated as well
as diffusion-dominated flows in an appropriate manner. We follow \cite{Lang1998} and set in each element
$K\in\cT_h$,
\begin{align}
\alpha_K &\;=\; \frac{h^\sharp_K}{2V} \frac{Re}{\sqrt{1+(Re)^2}},\quad
Re := \frac{h^\sharp_K V}{\mu},\\[1mm]
\beta_K &\;=\; \frac{h_K}{2\vert \bv \vert} \frac{Re}{\sqrt{1+(Re)^2}},\quad
Re := \frac{h_K \vert \bbeta \vert}{\mu},\\[1mm]
\beta_{K,i} &\;=\; \frac{h_K}{2\vert \bv \vert} \frac{Re}{\sqrt{1+(Re)^2}},\quad
Re := \frac{h_K \vert \bar{\bbeta} \vert}{\nu_i},\quad i=0,\ldots,n_s,
\end{align}
where $V$ denotes a global reference velocity, $h^\sharp_K$ is the diameter of the
two-dimensional ball having the same area as $K$, and $h_K$ ist the length of the
element $K$ in the direction of the local velocity $\bv$.

The highly nonlinear system (\ref{eq:var_form_fem}) is often hard to solve without
the knowledge of a good initial guess. A common technique is the use of a homotopy
approach in order to stabilize Newton's method. This can be best realized by a pseudo
implicit time-stepping scheme. We formally add time derivatives to equation
(\ref{eq:var_form_fem}), resulting in
\begin{equation}
\label{eq:var_form_fem_time}
\bu_h(t)\in V_h: \quad (P\partial_t \bu_h(t),\bphi) + A(\bu_h(t),\bphi) + S_h(\bu_h(t),\bphi)
= 0 \quad \text{for all } \bphi\in V_h
\end{equation}
with a projection matrix $P=\text{diag}(0,1,1)$. In this way, the
hydrodynamic pressure $p_{hyd}$ is
still determined by a stationary equation and therefore the characteristic property
of the low-Mach number model is preserved. Since an accurate resolution of $\bu_h(t)$ is not
important, we solve (\ref{eq:var_form_fem_time}) with the linearly implicit Rosenbrock
method {\sc Ros3pl} \cite{LangTeleaga2008} employing variable step sizes and a low tolerance
until a stationary solution is obtained. This method has very good stability properties and
is suitable for the solution of differential-algebraic equations like (\ref{eq:var_form_fem_time}).

\section{The Moving Mesh Method}
\subsection{Basic Formulation of the Moving Mesh PDE}
In what follows, we will adopt the time-dependent moving mesh method described
in \cite{HuangRenRussell1994,HuangRussell2010} to stationary combustion problems.
The mesh equation is formulated in terms of a coordinate transformation between
the original physical domain $\Omega$ and a computational domain $\Omega_C$ which
has the same topology. We denote the corresponding coordinates by
$\bx=(x_1,\dots,x_n)^T$ and $\bxi=(\xi_1,\dots,\xi_n)^T$. The time-depending mapping
$\bx(\bxi,t):\Omega_C\rightarrow\Omega$ is defined by the minimizer of
the quadratic functional
\begin{equation}
\label{eq:mm_functional}
I[\bxi] = \frac{1}{2} \int_{\Omega} \left( \sum_{i} \nabla \xi_i^T
\bG^{-1} \nabla \xi_i \right) d \bx,
\end{equation}
where $\boldsymbol G$ stands for the so-called monitor function that links
to the physical solution and is used to control the mesh density. Using a
variational approach, the minimizer is approximated by the smooth solution
of a modified gradient flow equation which reads
\begin{equation}
\label{eq:mmpde_xi}
\frac{\partial\bxi}{\partial t} = -\frac{B}{\tau} \frac{\delta I}{\delta \bxi}=
\frac{B}{\tau}\nabla \cdot ( \bG^{-1}  \nabla \bxi ).
\end{equation}
Here, $\tau>0$ is a user specified parameter that controls the response time
of the mesh movement and $B$ is a balance function to control the spatial movement
of the mesh points.
In order to compute the mapping $\bx(\bxi,t)$ directly, we interchange the roles
of dependent and independent variables. The final form of the moving mesh PDE (MMPDE)
can be expressed in the following form:
\begin{equation}
\label{eq:mmpde_x}
\tau  \frac{\partial \boldsymbol x}{\partial t}  = B \left( \sum_{i,j} A_{ij}
\frac{\partial^2 \boldsymbol x}{\partial \xi_i \partial \xi_j}
- \sum_{i} b_i \frac{\partial \boldsymbol x}{\partial \xi_i} \right),
\end{equation}
where
\begin{equation}
A_{ij} := \nabla \xi_i^T \boldsymbol G^{-1} \nabla \xi_j, \qquad \qquad
b_i:= \sum_j  \nabla \xi_i^T \frac{\partial \boldsymbol G^{-1}}{\partial \xi_j} \nabla \xi_j .
\end{equation}
The parameter $B$ should be determined such that all mesh points move with a
uniform time scale. This allows an easier numerical integration of the MMPDE.
The most convenient approach \cite{Huang2001} to obtain a well spatially balanced
MMPDE is to scale the terms on the right-hand side of (\ref{eq:mmpde_x}) as follows:
\begin{equation}
\label{eq:choice_of_P}
B = \frac{1}{\sqrt{\sum_i (A_{ii}^2 +b_i^2)}}.
\end{equation}
The main advantage of this approach is that the coefficients in the
right-hand side of (\ref{eq:mmpde_x}) become $\mathcal{O}(1)$. Then, the corresponding
equations are $\bx$- and $\bG$-scaling invariant, i.e., the MMPDE will not change if
we rescale the physical domain and the monitor function.  This allows to use the
positive parameter $\tau$ in (\ref{eq:mmpde_x}) to adjust the time scale of the
mesh movement. In general, a smaller $\tau$ results in a faster mesh adaptation
with respect to changes in the monitor function $\bG$, while a larger $\tau$ produces
slower movement in time \cite{CaoHuangRussell2001b,CaoHuangRussell2001a}. Since we
only focus on the design of a static quasi-optimal mesh to improve the approximation
quality of the stationary solution, it is sufficient to set $\tau=1$ and to apply a
pseudo time-stepping method. For a complete specification of the coordinate transformation,
we also need to supply the MMPDE with appropriate boundary conditions. In our applications,
we always fix the mesh points on the boundary. For more advanced strategies,
we refer to \cite{Huang2001}.

\subsection{Construction of the Monitor Function}
A general approach for constructing monitor functions $\bG$ in two
dimensions is based on its eigendecomposition,
\begin{equation}
\label{eq:def_G}
\bG = \lambda_1\bitv_1\bitv_1^T + \lambda_2\bitv_2\bitv_2^T,
\end{equation}
where $\bitv_1$ and $\bitv_2$ are mutually orthogonal normalized
eigenvectors of $\bG$. The eigenvectors of $\bG$ determine the
directions of the mesh adaptation, while the associated eigenvalues dictate
the intensity of the concentration of the mesh in these directions.

Given a vector-values quantity of interest $\bpsi(\bu_h(\bx))$ computed from
the stationary numerical solution $\bu_h(\bx))$, a class of monitor functions
can be constructed by harmonic mappings \cite{CaoHuangRussell1999}. This reads
\begin{equation}
\label{eq:def_lambdaForG}
\begin{array}{l}
\ds\bitv_1 = \frac{\bpsi}{\Vert\bpsi\Vert}, \quad \bitv_2 = \bitv_1^{\bot},\\[5mm]
\ds\lambda_1 = \sqrt{1+\alpha\Vert\bpsi\Vert^2},
\qquad \lambda_2 = \frac{1}{\lambda_1},
\end{array}
\end{equation}
where $\alpha$ denotes a user-defined intensity parameter. Typical examples for
the function $\bpsi$ are error indicators, curvature or gradient data taken from
certain solution components. Heuristically, the numerical error is larger in regions
where the solution changes dramatically. For example, suppose a numerical scalar solution
$u_h$ exhibits a steep front and the solution gradient is chosen to be one of the
eigenvectors, i.e., $\bpsi=\nabla u_h$, then, by using (\ref{eq:def_lambdaForG}),
the intensity of the adaptation determined by $\lambda_1$ along the normalized gradient
direction $\bitv_1$ is much stronger than those determined by $\lambda_2$ in
the tangential direction $\bitv_2$. Thus, it is expected that coordinate expansion
and compression will mainly occur in the gradient direction. However, it should be
pointed out that using the solution gradient as $\bpsi$ may not always be the best
option for some problems. This topic will be discussed in our
applications. Similar conclusion can already be found in \cite{CaoHuangRussell2001b}.

The monitor function $\bG$ is usually constructed by non-smooth operations applied
to the numerical solutions. This can result in very stiff MMPDEs. A usual remedy
is to smooth the monitor function. Let $\bx_p$ be a mesh point in $\Omega$ and $\bxi_p$
the corresponding mesh point in $\Omega_c$. Then the following smoothing algorithm has
proven to work quite satisfactorily in practice \cite{Huang2001,LangCaoHuangRussel2003}:
\begin{equation}
\begin{array}{ll}
\bG^0(\bx_p) :=\bG(\bx_p) & \text{ for all grid points } \bx_p,\\[0.7em]
\ds\bG^{m+1}(\bx_p) := \frac{\int_{C(\bxi_p)} \bG^m(\bx_p(\bxi))\,d\bxi}{\vert C(\bxi_p)\vert }
& \text{ for all grid points } \bx_p, \\
& \;m=0,1,\ldots, M_s-1,
\end{array}
\end{equation}
with $C(\bxi_p) \subset \Omega_C$ is the union of neighbouring grid cells having
$\bxi_p$ as their common vertex. Here, the initial monitor function $\bG(\bx_p)$
is directly obtained from (\ref{eq:def_G}) and (\ref{eq:def_lambdaForG}).
$M_s$ is a user-specified parameter standing for the number of the smoothing cycles to
be performed. A smaller $M_s$ gives a more accurate description of the profiles of the
monitor function, since it characterizes the exact features of the solutions. However,
this usually makes the MMPDE harder to solve.  Generally, the value of the intensity
parameter $\alpha$ and $M_s$ are problem-dependent.

\section{Applications}
In this section, we present two reactive flow problems which have an increasing degree of complexity. Both problems have been extensively studied in \cite{Braack1998,BraackRannacher1999}. The first example is the ozone decomposition flame. The chemical processes are modelled with $3$ species and $6$ elementary chemical reactions. The motivation for this application is to demonstrate and compare the numerical results for various monitor functions. Benefits and shortcomings of the moving mesh method will also be discussed. The second example describes a methane flame in a complex geometry. The underlying chemical model contains $15$ species and $84$ elementary reactions. Some intermediate species have extremely thin flame layers that need to be resolved sufficiently. Furthermore, the whole process possesses both extremely fast and slow motions, implying a strict requirement on the time resolution.
All transport coefficients for viscosity, thermal conductivity and diffusion are evaluated
from kinetic models collected in the data bases of the Sandia National Laboratories \cite{KeeRupleyMiller1987}. The influence of the gravitation is neglected, i.e., we
set $\bg=0$ in our applications.

\subsection{A Two-Dimensional Example for Ozone Decomposition}
We use the mechanism for ozone decomposition reaction  that consists of $6$ elementary reactions proposed in \cite{Warnatz1978}:
\begin{eqnarray}
\label{eq:Reaction_O3_1}
\text{O}_3 +  \text{M}
&  \overset{1}{\underset{2}{ \rightleftharpoons }}&    \text{O}_2 + \text{O} + \text{M}\\
\text{O} + \text{O} + \text{M} & \overset{3}{\underset{4}{ \rightleftharpoons }}&  \text{O}_2  + \text{M}\\
\label{eq:Reaction_O3_3}
\text{O}+ \text{O}_3  & \overset{5}{\underset{6}{ \rightleftharpoons }}&  \text{O}_2 + \text{O}_2
\end{eqnarray}
Here, $\text{M}$ denotes an arbitrary third body, i.e., it can be one of three considered species: Oxygen atoms $\text{O}$, Oxygen molecules $\text{O}_2$ or Ozone $\text{O}_3$.
Further details of the reaction mechanism based on the Arrhenius law can be found in
\cite[Tab. A.1]{BraackRannacher1999}.

\subsubsection{Specification of the Simulation}
The geometry of the simulation is supposed to consist of two flat plates (regarded as rigid walls) with an inflow at the left side, comprising a cold mixture of ozone and oxygen molecules, and a free-stream outflow at the right side. The computational domain is defined as $\Omega:=]0, 0.02[ \times ]0, 0.005[$ and $\partial \Omega=\Gamma_{\text{in}}   \cup \Gamma_{\text{wall}}  \cup \Gamma_{\text{out}}$. The boundary conditions are
\begin{eqnarray}
\text{inflow on }\Gamma_{\text{in}}\text{:}
&& \bv =\bv_{in}, \,  T=T_{in},\  \bw =\bw_{in}, \\
\text{rigid walls on } \Gamma_{\text{wall}}\text{:}
&& \bv = 0, \,  T=T_{wall}, \, \nabla\bw\cdot\bn = 0,  \\
\text{outflow on }\Gamma_{\text{out}}\text{:} &&
-\mu \nabla\bv\cdot\bn + p \bn = 0, \, \nabla T\cdot\bn = 0,
\nabla\bw\cdot\bn = 0.
\end{eqnarray}
The initial mass fractions of species are $\omega_{\text{O}_3}=0.2$, $\omega_{\text{O}_2}=0.8$
and the initial temperature is
\begin{equation}
T(x) = 298 + (800-298)\,e^{-10^5(x_1-0.005)^2},
\end{equation}
which is kept fixed at the inflow and at the rigid wall as Dirichlet boundary conditions.
The velocity profile on the inflow boundary is parabolic with a maximum velocity of $0.25$  m/s.  The Reynolds number is approximately $Re=v_{c}L/\mu\approx 62$, where $v_c$ is the mean value of the velocity at the inlet, $L$ stands for the length of the geometry, and ${\mu}$ is the viscosity of the mixture at the inflow.

\begin{figure}[t]
\centering
\includegraphics[width=0.7\linewidth]{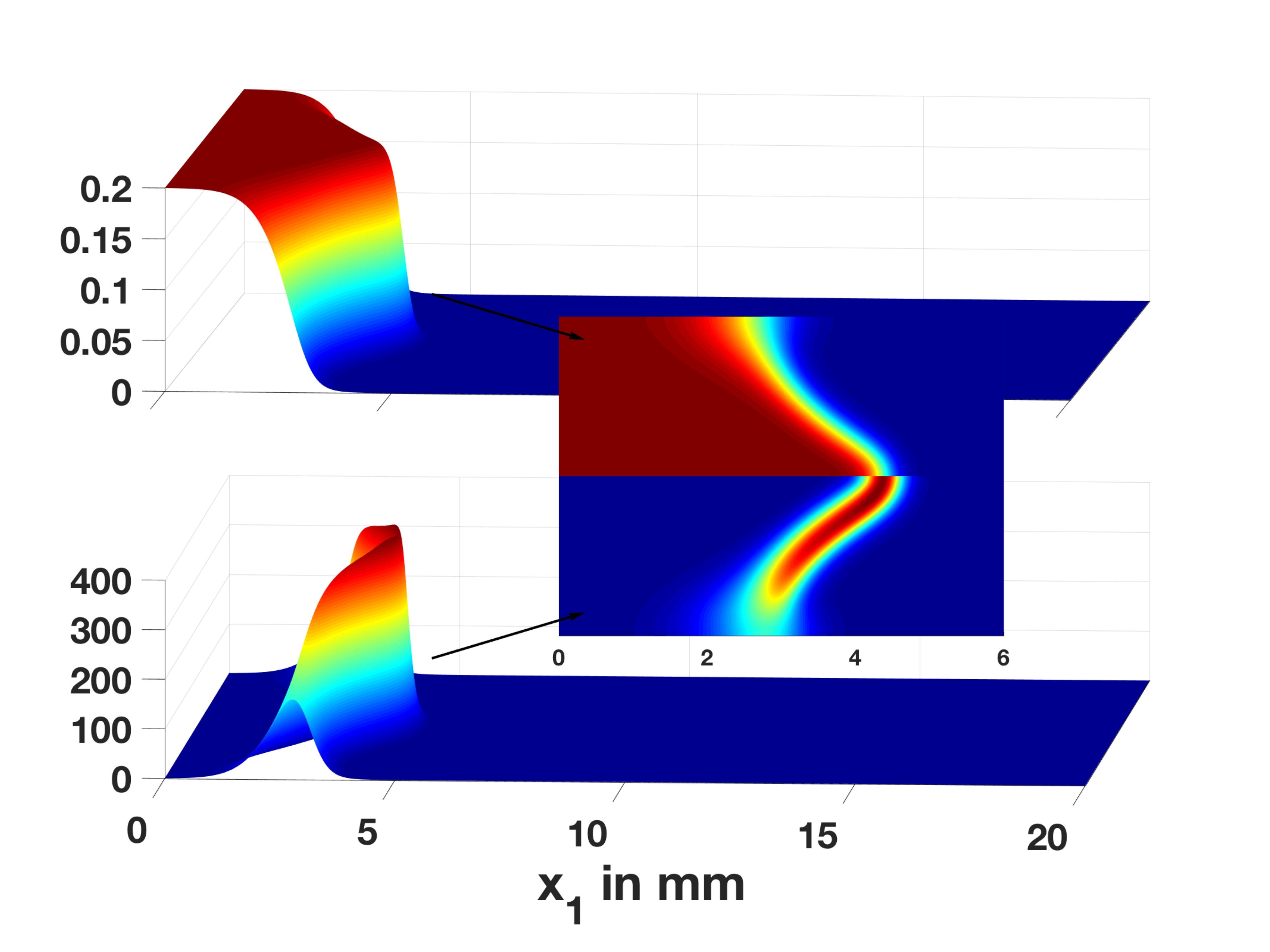}
\parbox{13cm}{
\caption{Profiles of the numerical solution $w_{\text{O}_3,h}$ (top)
and its recovered gradient $\vert D^1 w_{\text{O}_3,h}\vert$ (bottom).}
\label{fig:FirstGradientOfW03Demo}
}
\end{figure}

A good indicator of the correct location of the flame front is the
mean value of $\omega_{\text{O}_3}$ \cite{BraackRannacher1999},
\begin{equation}
\label{qoi:TargetOzone}
J(\bu)= \frac{1}{\vert \Omega \vert}\int_{\Omega} w_{\text{O}_3} d\bx.
\end{equation}
The profiles of the numerical approximation of $w_{\text{O}_3}$ and its
gradient are shown in Fig. \ref{fig:FirstGradientOfW03Demo}.
We will use $J(\bu_h)$ to demonstrate the performance of our moving
mesh approach. To compare the accuracies for different meshes, we first
compute a reference solution with an adaptive, locally refined mesh with
about $45000$ grid points using the {\sc Kardos} software
\cite{ErdmannLangRoitzsch2002,Lang2000} and provide also values for
regular meshes with different numbers of uniform grid points. The
results are collected in Tab.~\ref{tab:numerical_Result_w_o3}.
\begin{table}[h!]
\centering
\begin{tabular}{ r | r | r }
\hline
noP  & $J(\bu_h) $  &  $ \vert J(\be) \vert$\\
\hline\rule{0mm}{5mm}
1105 & 0.03331313 & 1.912892 $\times 10^{-4}$ \\
2835 & 0.03347317 & 1.570370 $\times 10^{-4}$ \\
4257 & 0.03348598 & 1.844070 $\times 10^{-5}$ \\
\hline\rule{0mm}{5mm}
adaptive & 0.03350442 & -\\
\hline
\end{tabular}\\[2mm]
\parbox{13cm}{
\caption{Values for $J(\bu_h)$ defined in (\ref{qoi:TargetOzone}) for uniform
meshes with different numbers of points (noP) and absolute
values of the numerical errors $J(\be)=J(\bu-\bu_h)$ derived
from an adaptive reference solution with about $45000$ grid points.  }
\label{tab:numerical_Result_w_o3}
}
\end{table}

\subsubsection{A First Monitor Function}
\begin{figure}[htp]
\centering
\includegraphics[width=0.46\linewidth]{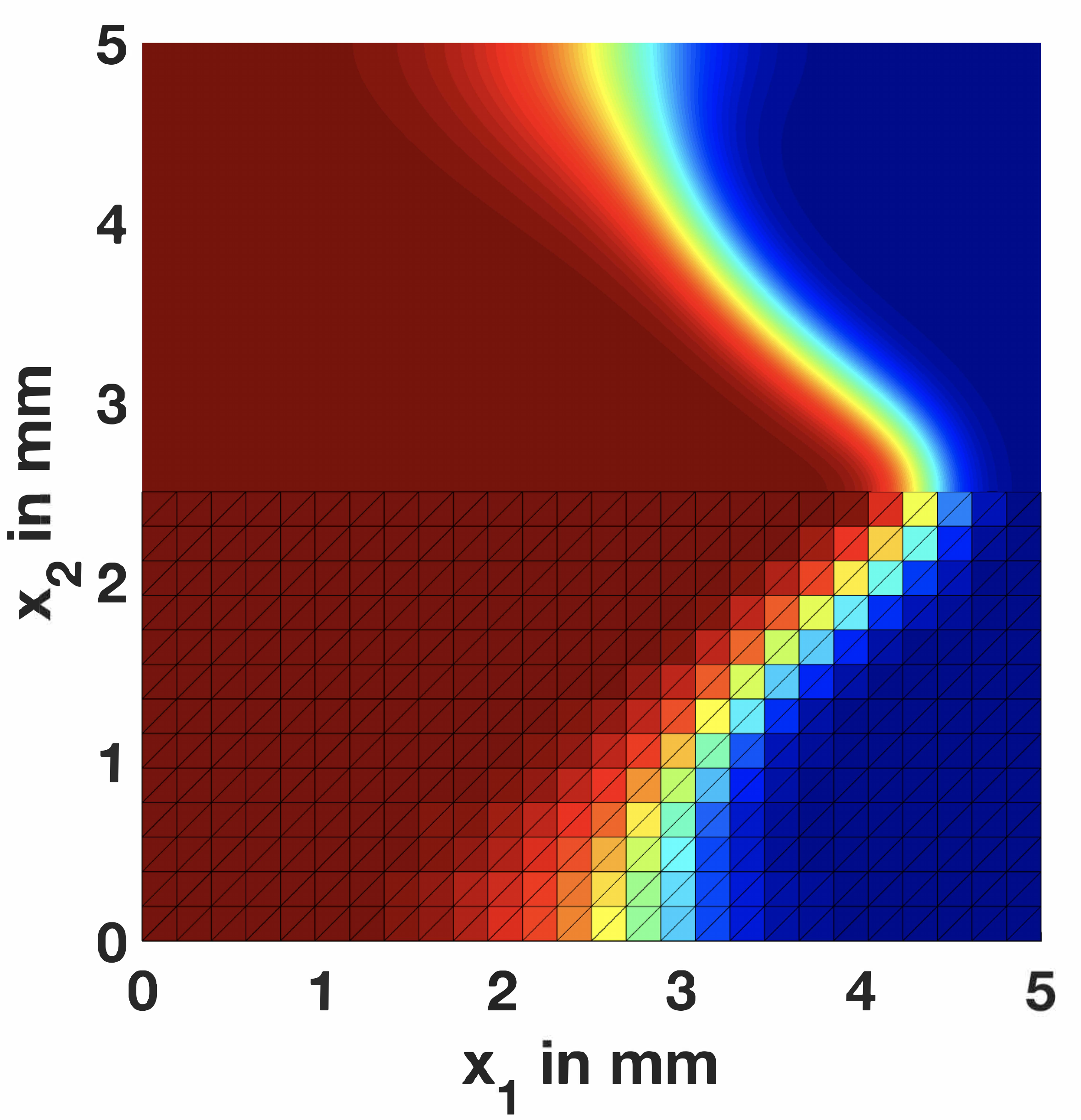}
\includegraphics[width=0.46\linewidth]{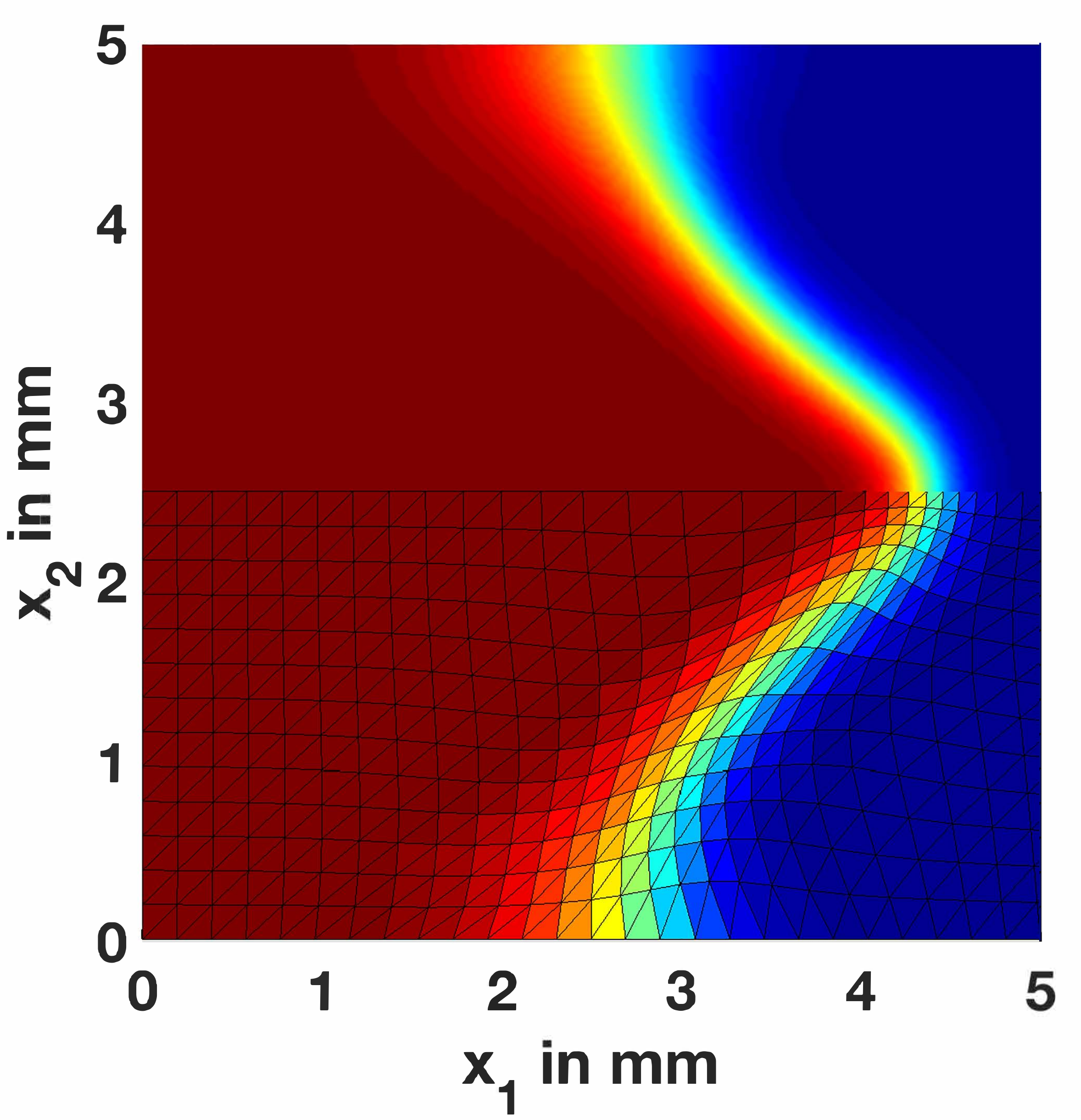}\\[2mm]
\parbox{13cm}{
\caption{Comparison of the numerical solution $w_{\text{O}_3,h}$ computed on uniform
and adaptive meshes. Left: numerical solution $w_{\text{O}_3,h}$ with high resolution computed on a uniform mesh with 148417 grid points (above) and with low resolution on a uniform mesh with 2835 grid points (below). Right: numerical solution $w_{\text{O}_3,h}$ with high resolution computed on a regular mesh with 148417 grid points (above) and on an adaptive mesh with 2835 grid points (below).}
\label{fig:ozone-m1-wo3}
}
\end{figure}

\begin{figure}[htp]
\centering
\includegraphics[width=0.46\linewidth]{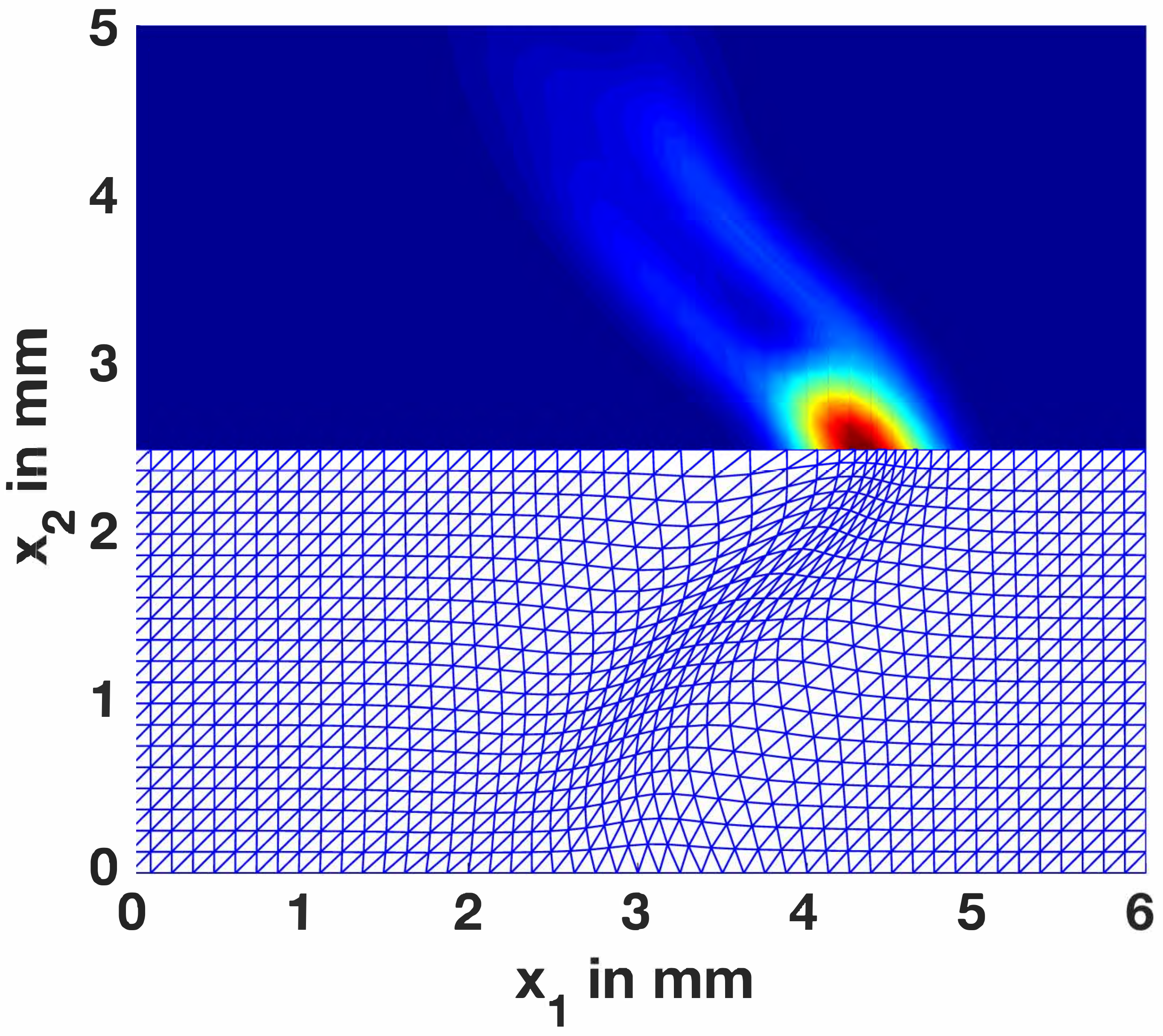}
\includegraphics[width=0.50\linewidth]{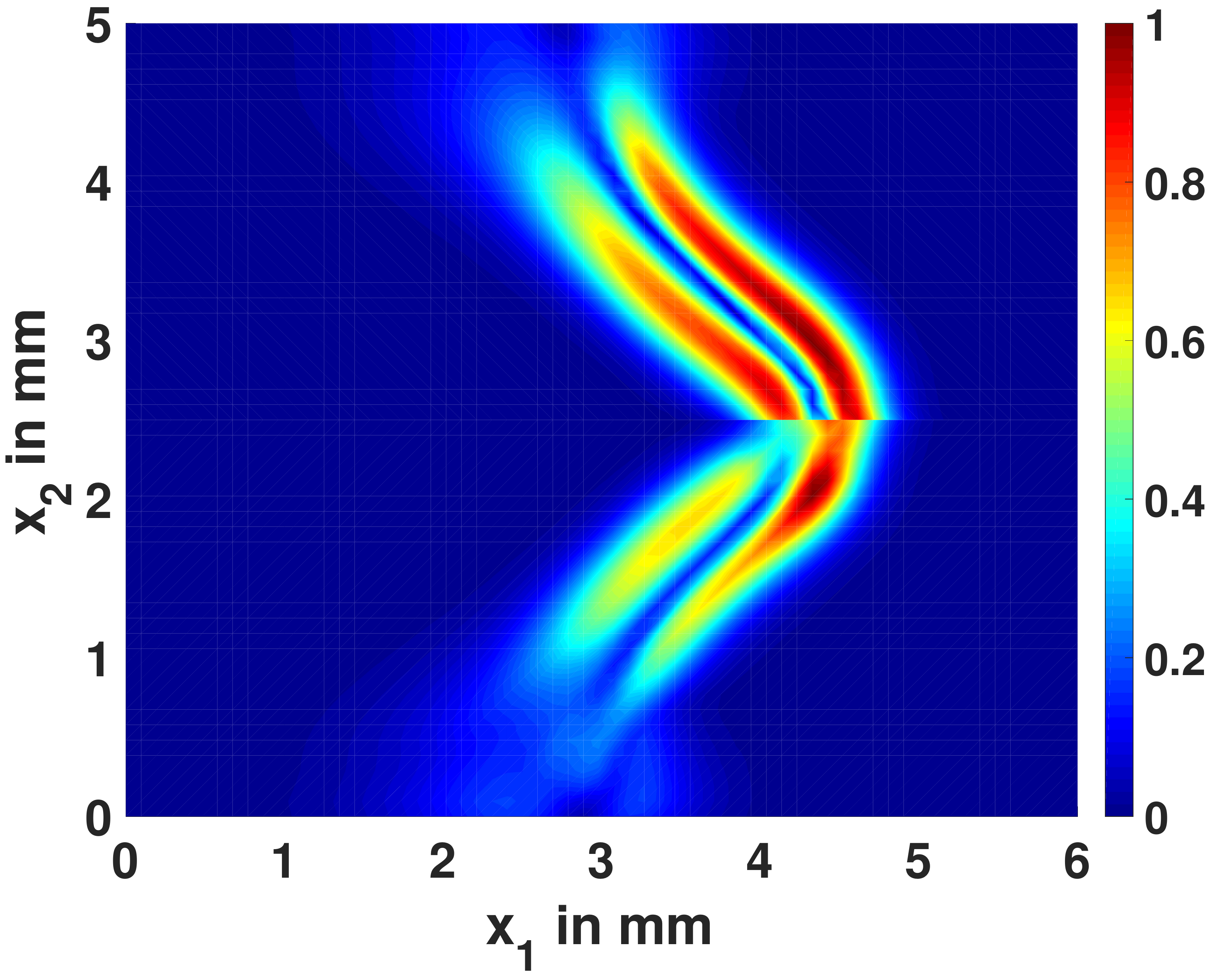}\\[2mm]
\parbox{13cm}{
\caption{Left: 2D plot of the mesh density $\sqrt{\det (\bG(D^1 w_{\text{O}_3,h}(\bx)))}$ after
smoothing (above) and corresponding adaptive moving mesh (below). Right: error distribution of the mass fraction of ozone molecules $\omega_{\text{O}_3}$; approximate interpolation error indicator (above) and hierarchical basis error estimator (below). }
\label{fig:ozone-m1-wo3-density}
}
\end{figure}

A central issue when applying a moving mesh method is to define a proper monitor function to control the distribution of the grid points. In most cases, the monitor function should ensure that the grid points are concentrated in regions where the physical solutions require a higher resolution. So the grid points on a uniform mesh should be relocated so that the requirement on the number of the mesh cells necessary to appropriately resolve a flame layer can be satisfied \cite{PoinsotVeynante2012}.  Common practice when constructing a monitor function is to employ some heuristic choices. For problems where the solutions change dramatically within a small region, such as the flame layer of the ozone molecules $\omega_{\text{O}_3}$ illustrated in Fig. \ref{fig:FirstGradientOfW03Demo}, there are usually larger numerical errors in regions with large gradients. Thus the gradients of certain solution components often represent a good error indicator. The gradient of a piecewise linear function is a constant vector in each
element $K\in\cT_h$. In order to improve the approximation property, we apply a standard
gradient recovery operator $D^1: S_h\rightarrow (S_h)^2$ as introduced in \cite{ZienkiewiczZhu1992}, and choose in our first experiment
\begin{equation}
\bpsi : = D^1 w_{\text{O}_3,h}
\end{equation}
as eigenvector $\bv_1$ of the monitor function $\bG$ to enforce a mesh adaption
in this direction. From Fig. \ref{fig:ozone-m1-wo3}, we can make the following observations: The contour of the flame layer becomes more distinct and can be better approximated through adjusting the locations of the grid points close to the sharp flame front.  The mesh cells are compressed to different extents along the convective direction of the flow so that an obvious anisotropic behaviour is visible.  Moreover, the mesh points are correctly concentrated in the area where the gradient is large. This is exactly the location where the mesh density takes its peak values, as illustrated in Fig. \ref{fig:ozone-m1-wo3-density}. However,
Tab.~\ref{tab:ozone-m1-j} shows that the accuracy of $J(\bu_h)$ does not always improve with increasing grid points as one would expect. A similar behaviour has been observed in \cite{Huang2001}.

\begin{table}[h!]
\centering
\begin{tabular}{|r|r|r|r|r|}
\hline
noP & $J(\bu_h)$, uniform & $J(\bu_h)$, adaptive
& $\vert J(\be) \vert$, uniform &  $\vert J(\be) \vert$, adaptive \\
\hline\rule{0mm}{2mm}
1105 & 0.03331313 & 0.03356365 & 2.912892 $\times 10^{-4}$ & 5.918997 $\times 10^{-5}$ \\
\hline\rule{0mm}{2mm}
2835 & 0.03347317 & 0.03344742 & 1.570370 $\times 10^{-4}$ & 5.703389 $\times 10^{-5}$ \\
\hline\rule{0mm}{2mm}
4257 & 0.03348598 & 0.03344386 & 1.844070 $\times 10^{-5}$ & 6.059835 $\times 10^{-5}$ \\
\hline
\end{tabular}\\[2mm]
\parbox{13cm}{
\caption{Results for uniform and moving meshes with the quantity of
interest $\bpsi=D^1 w_{\text{O}_3,h}$ and parameters $M_s=8$ and
$\alpha=80$ for the moving mesh design. The reference value is
$J(\bu)\approx 0.03350442$.}
\label{tab:ozone-m1-j}
}
\end{table}
In the following, we would like to discuss two reasons for this behaviour. First,
while the grid points are concentrated in regions with large gradients of
$w_{\text{O}_3,h}$, the accuracy of the other components, which also influence the
balance of $\text{O}_3$, might be lower due to certain local expansions of mesh lines,
resulting in a worse accuracy for $J(\bu_h)$. Second, although regions of large solution gradients can usually reflect areas with large numerical errors, this is not
always true for the exact error distribution. Because of the coupling of all components, using the gradient of a single component may lead to an over-concentration of mesh points in areas that have already been sufficiently resolved. The right graph in
Fig.~\ref{fig:ozone-m1-wo3-density} illustrates an estimated error distribution which is based on a hierarchical error estimator proposed in \cite{Lang2000}. Obviously, the error distribution differs significantly from the density function illustrated in the left graph in Fig.~\ref{fig:ozone-m1-wo3-density}, which governs the concentration of the grid points in the MMPDE. The error distribution has a double-layer structure, where the numerical errors are mainly accumulated in regions with large flame curvature.  More specifically, at the downstream edge, the numerical errors are larger than those on the upstream side. The downstream edge of the flame layer is located in a transitional reaction zone, where the dynamics of the species is governed by both convection and source terms. There the intermediate species are also frequently produced and consumed. Thus, a higher resolution in this region is preferable.
Next we will construct a monitor function which is based on local error information.

\subsubsection{A Second Monitor Function}
\begin{figure}[htp]
\centering
\includegraphics[width=0.48\linewidth]{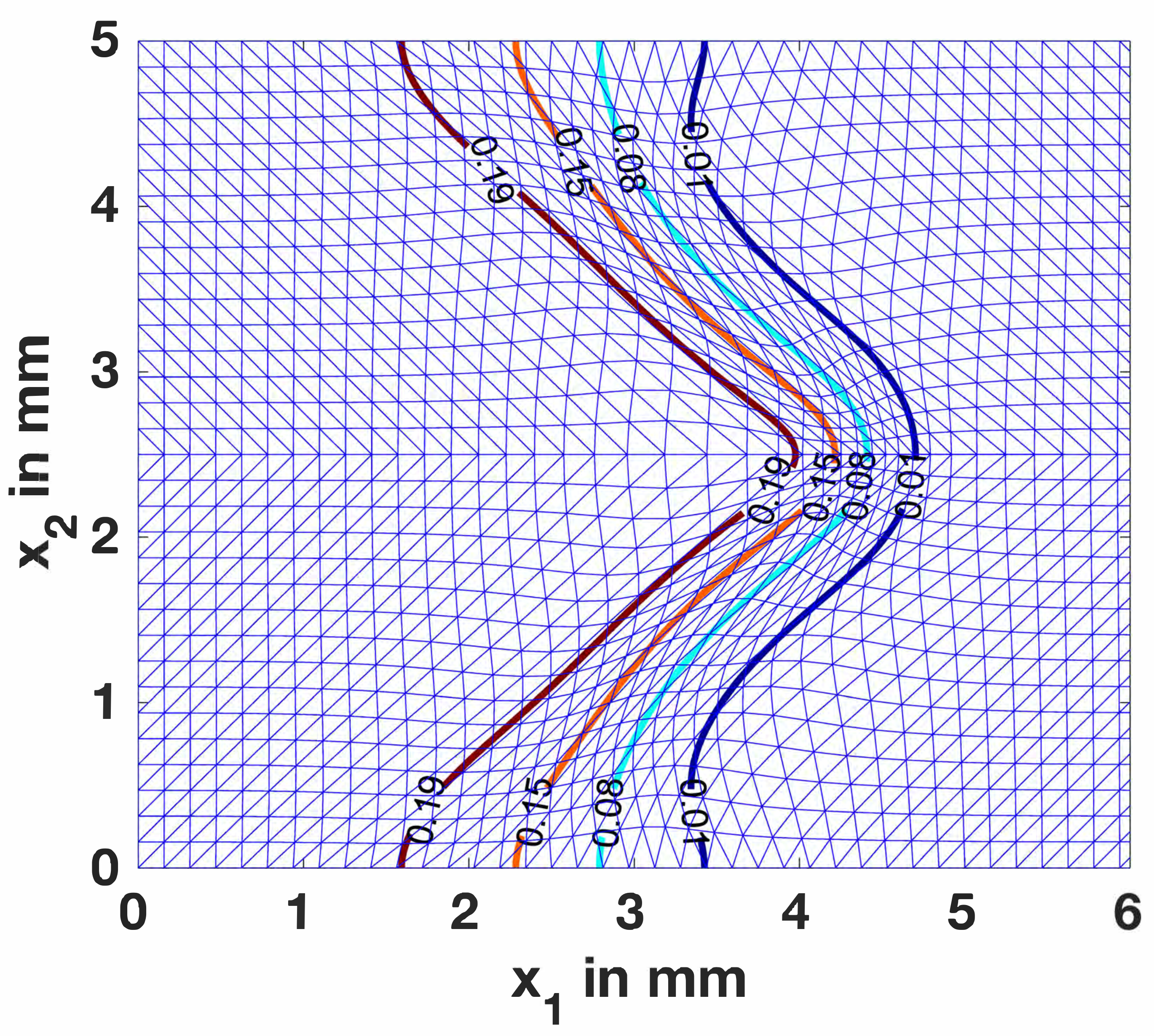}
\includegraphics[width=0.48\linewidth]{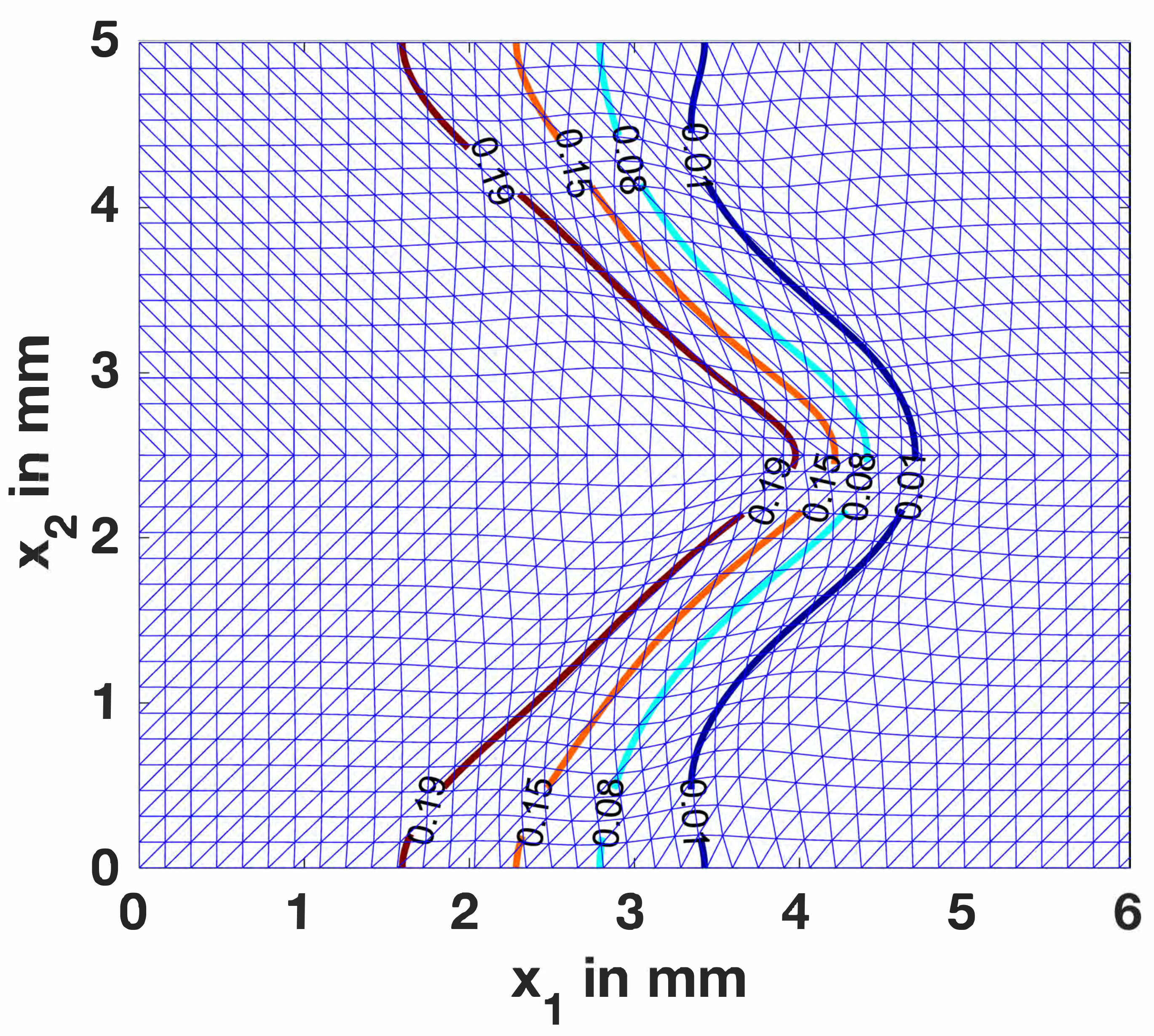}\\[2mm]
\parbox{13cm}{
\caption{Comparisons of moving meshes generated by different monitor functions
for $4257$ mesh points. Colour curves represent the contour lines of $\omega_{\text{O}_3}$. Left: Close-up view of the adaptive mesh around the O$_3$-flame layer computed by the gradient-based function $\bpsi=D^1 w_{\text{O}_3,h}$.  Right:
Close-up view of the adaptive mesh around the O$_3$-flame layer computed by the
second derivative-based function $\bpsi=D^2 w_{\text{O}_3,h}$.}
\label{fig:ozone-m1m2-wo3-comparison}
}
\end{figure}

Now, we will introduce a second monitor function which is based on an approximate
interpolation error indicator. First, we recall the Sobolev space $H^2(\Omega)$ which
consists of all functions from $H^1(\Omega)$ with square-integrable generalized second-order derivatives. In what follows, we will denote by $|\cdot|_{H^i}, i=1,2$, the usual semi-norm in
$H^i$ defined by derivatives of $i$-th order only.
Let $u\in H^2(\Omega)\cap C^0(\bar{\Omega})$, for the moment, be a scalar solution,
$u_h$ its linear finite element approximation, and $\mathnormal{\Pi}_1u\in S_h$
its linear interpolant defined by $u(\bx_i)=\mathnormal{\Pi}_1u(\bx_i)$ for all
mesh points  $\bx_i$. Then a standard argument for the interpolation error
\begin{equation}
\label{eq:interpolation}
\vert u - \mathnormal{\Pi}_1 u \vert_{H^1(\Omega) } =
\left( \sum_{K \in \cT_{h}} \vert u - \mathnormal{\Pi}_1 u \vert_{H^1(K)}
\right)^{1/2} \leq  C \sum_{K \in \mathcal{T}_{h}} h_K \vert u \vert_{H^2(K)},
\end{equation}
gives the a priori error estimate for the linear finite element solution,
\begin{equation}
\label{eq:h1error}
\Vert u - u_h \Vert_{H^1(\Omega) } \leq  C h \vert u \vert_{H^{2}(\Omega)},
\end{equation}
where $h:=\max_{K\in\cT_h}h_K$. This, together with the observations made in
the previous section, motivates to use second derivatives of $u_h$ to control
the local density of the moving mesh. An approximation $D^2u_h\in (S_h)^2$ for
$(\partial_{x_1x_1}u,\partial_{x_2x_2}u)^T$ is derived from a quadratic interpolation
of $u_h$ \cite[Sec. 4.2]{CaoHuangRussell2001b}. Then, we set for our application,
\begin{equation}
\bpsi : = D^2 w_{\text{O}_3,h}
\end{equation}
\begin{figure}[t]
\centering
\includegraphics[width=0.8\linewidth]{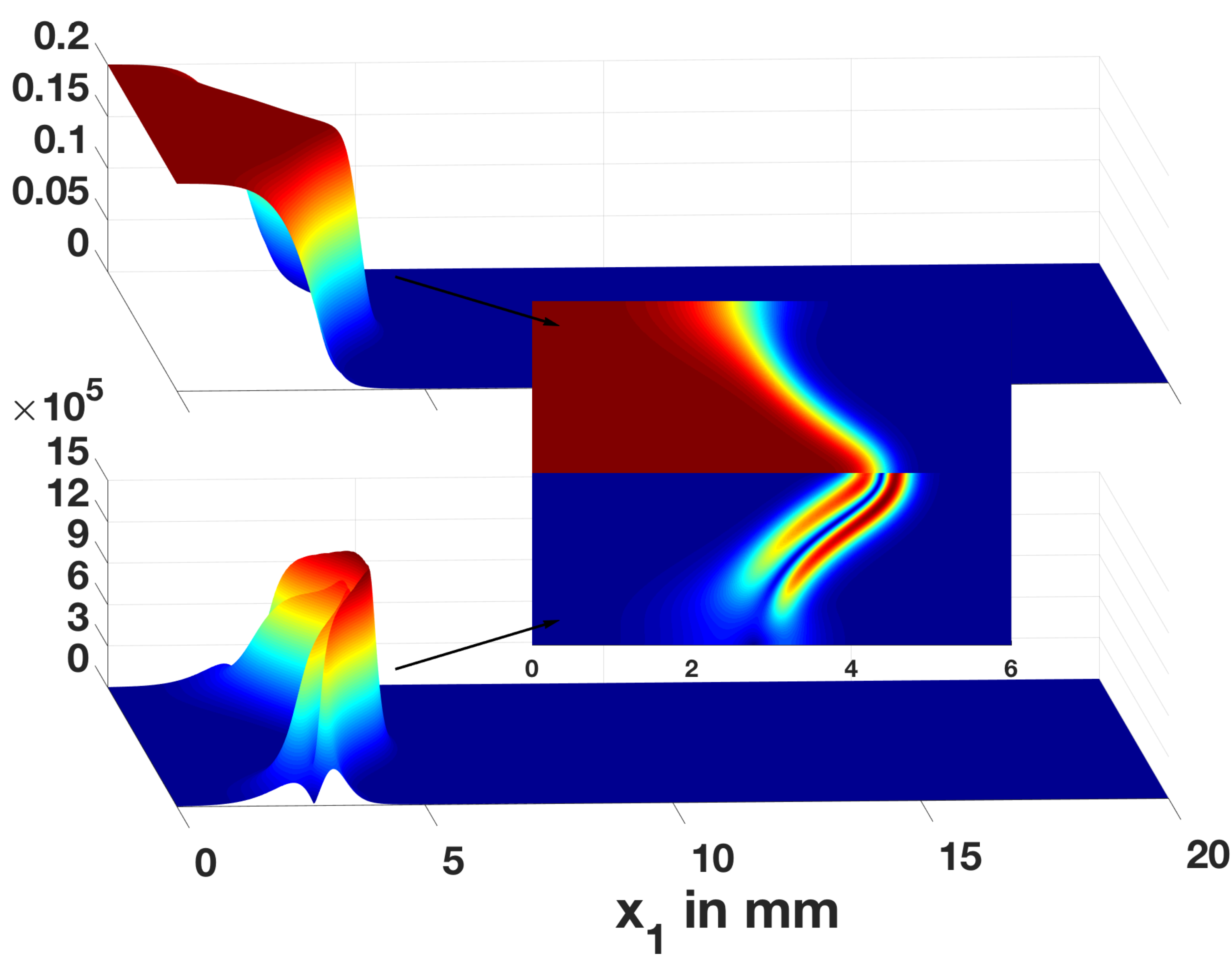}
\parbox{13cm}{
\caption{Profiles of the numerical solution $w_{\text{O}_{3,h}}$ (top)
and its interpolation error indicator $\vert D^2 w_{\text{O}_{3,h}}\vert$ (bottom).}
\label{fig:FGradientOfW03Demo}
}
\end{figure}
\begin{table}[h!]
\centering
\begin{tabular}{|r|r|r|r|r|}
\hline
noP  & $J(\bu_h)$, gradient &$J(\bu_h)$, interr &
$\vert J(\be)\vert$, gradient & $ \vert J(\be) \vert$, interr \\
\hline\rule{0mm}{2mm}
2835 & 0.03344742 & 0.03344717 & 5.703389 $\times 10^{-5}$ & 3.272824 $\times 10^{-5}$ \\
\hline\rule{0mm}{2mm}
4257 & 0.03344386 & 0.03504592 & 6.059835 $\times 10^{-5}$ & 1.347509 $\times 10^{-7}$ \\
\hline
\end{tabular}\\[2mm]
\parbox{13cm}{
\caption{Computation results for moving meshes with two monitor function based on the gradient
$D^1 w_{\text{O}_{3,h}}$ (gradient) and on the interpolation error indicator $D^2 w_{\text{O}_{3,h}}$ (interr) with $M_s=8$,
$\alpha=80$. The reference value is $J(\bu)\approx 0.03350442$.}
\label{tab:ozone-m1m2-j}
}
\end{table}
In Fig.~\ref{fig:FGradientOfW03Demo}, the profile of $w_{\text{O}_3,h}$ and its
approximate interpolation error indicator $\vert D^2 w_{\text{O}_3,h} \vert$
are shown. We observe that this error indicators gives nearly as good information
on the local errors as the more sophisticated error estimator based on a hierarchical
basis, compare Fig.~\ref{fig:ozone-m1-wo3-density}, but is much cheaper. From Tab.~\ref{tab:ozone-m1m2-j}, we see that the results obtained by using the new
monitor function are much better than those based on the gradient of $w_{\text{O}_3,h}$.
From a comparison of both methods shown in Fig~\ref{fig:ozone-m1m2-wo3-comparison},
it becomes visible that with the new monitor function grid points are concentrated in
a larger area due to the higher second derivatives and hence numerical errors.

\subsection{A 2D Example for Methane Combustion}
In this section, we consider a methane combustion problem described by
the following global reaction
\begin{equation}
\text{CH}_4 + 2 \ \text{O}_2  \rightarrow \text{CO}_2  + 2  \ \text{H}_2 \text{O}.
\end{equation}
The reaction mechanism with $15$ species and $84$ elementary reactions is
taken from \cite{Smooke1991}. We will consider a lamella burner which is modified from a
prototype constructed by JUNKERS Bosch Thermotechnik \cite{ParmentierBraackRiedelWarnatz2003}.
The geometry for the simulation is given in Fig.~\ref{fig:GeometryForBurnerofCombustionOfMethane}.

\subsubsection{Specification of the Simulation}
\begin{figure}[t]
\centering
\includegraphics[width=0.7\linewidth]{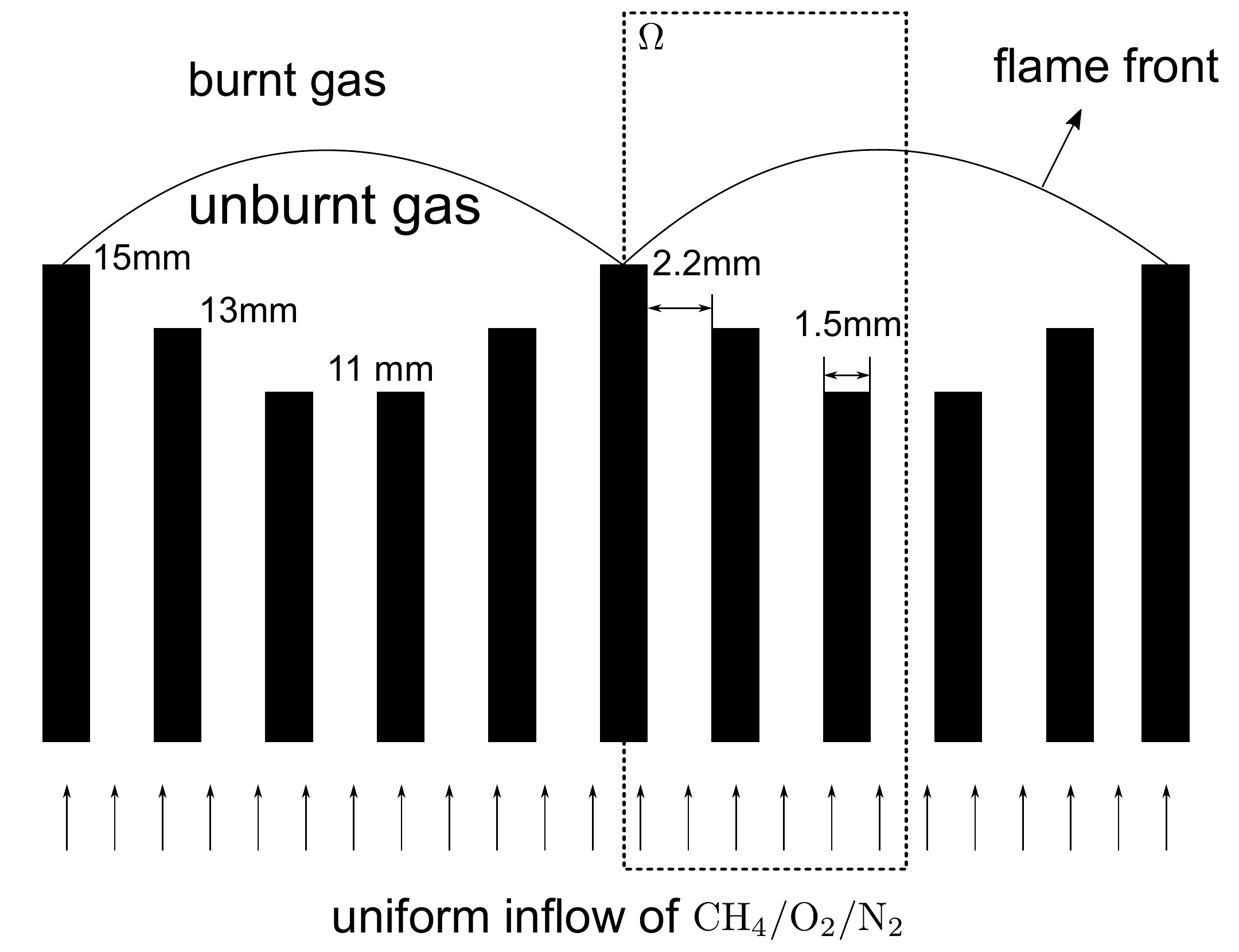}\\[2mm]
\parbox{13cm}{
\caption{Geometry of the methane burner considered.}
\label{fig:GeometryForBurnerofCombustionOfMethane}
}
\end{figure}

A stoichiometric mixture of methane (CH$_4$) and air (O$_2 /$N$_2$) flows from the bottom of
the burner through a sample of slots which have a uniform width interval of $2$ mm and three different heights varying from $15$ mm to $11$ mm. All slots have the same width of $1.5$ mm. The inflow velocity is $v_2=0.27$ m/s. The solution is assumed to be spatially periodic. Thus it is sufficient to restrict the computational domain to three slots which define $\Omega$, see
Fig.~\ref{fig:GeometryForBurnerofCombustionOfMethane}. The lamellae can be considered as obstacles. Dirichlet boundary conditions are used for the temperature, no-slip conditions for the velocities and Neumann boundary conditions for the species on the wall of the lamellae. To specify the temperature on the wall of the slots, we use three piecewise linear functions varying from $298$ K to $393$ K, $453$ K, $463$ K, respectively. On the cut boundary of $\Omega$, symmetric boundary conditions are used for all the variables. For more details, we
refer to \cite{Braack1998,BraackRannacher1999,Sun2018}.

This problem is characterized by the interaction of different physical processes and
largely separated scales. A simulation on uniform meshes can be very prohibitive.
We expect that moving mesh strategies designed with an appropriate monitor function will improve the numerical approximation even on relatively coarse meshes.

\subsubsection{Application of Moving Meshes}
Since the underlying combustion mechanism is relatively complicated, it is
not possible to calculate a stable stationary solution on quasi-uniform
meshes with around $40000$ grid points. Higher resolution is needed near the
flame front in the reaction zone and in the pre-heating zone, where the chemical
and convection-diffusion processes are strongly coupled. These regions
are unknown a priori. Therefore, we directly couple the physical PDE with the
MMPDE by an Arbitrary Lagrangian-Eulerian approach.

We first define $\bhu_h(t):=\bhu_h(\bx(\bxi,t),t)=\bhu_h(\bxi,t)$. Under the mapping
$\bx(\bxi,t)$, we then transform (\ref{eq:var_form_fem_time}) into
a system involving the computational coordinates $\bxi$, i.e.,
\begin{equation}
\label{eq:var_form_fem_time_cdom}
\begin{array}{rll}
\bhu_h(t)\in V_h: && \left( P \left( \partial_t \bhu_h(t) -
(\partial_t\bx\cdot J^{-T} \hat{\nabla} )\,\bhu_h(t)\right),\bphi \right)\\[3mm]
&& + \;\hat{A}(\bhu_h(t),\bphi) + \hat{S}_h(\bhu_h(t),\bphi) = 0 \quad\quad
\text{for all } \bphi\in V_h.
\end{array}
\end{equation}
Here, $\hat{\nabla}$ denotes the gradient operator with respect to $\bxi$ and
$J=\partial\bx/\partial\bxi$ is the Jacobian of the mapping $\bx(\bxi,t)$.
The operators $\hat{A}$ and $\hat{S}_h$ are obtained from $A$ and $S_h$ by
using the identity $\nabla = J^{-T}\hat{\nabla}$ and replacing $\bu_h$ by
$\bhu_h$.
The additional term $(\partial_t\bx\cdot J^{-T} \hat{\nabla})\,\bhu_h$ on the left
side can be viewed as a correction for the convective effects of the mesh motion.
Note that the linear finite element space $V_h$ is now related to a
time-independent regular
partition of triangles of the computational domain $\Omega_C$. Equation
(\ref{eq:var_form_fem_time_cdom}) is simultaneously solved with the MMPDE defined
in (\ref{eq:mmpde_x}) in the course of the pseudo time-stepping method. The location
of the mesh points described by the mapping $\bx(\bxi,t)$ is
immediately adapted to the evolving flame structure until the stationary state is
reached. Adjusting the mesh points in this way yields a stable numerical solution,
even on relatively coarse meshes. In our calculation, we start with an isotropic
quasi-uniform mesh with 29864 grid points delivered from the 2D-mesh
generator {\sc Triangle} \cite{Shewchuk2002}. We first compute the
temperature-dependent flow field without combustion and start then our mesh moving
approach including all reaction terms.

\begin{figure}[ht]
\centering
\includegraphics[width=0.9\linewidth]{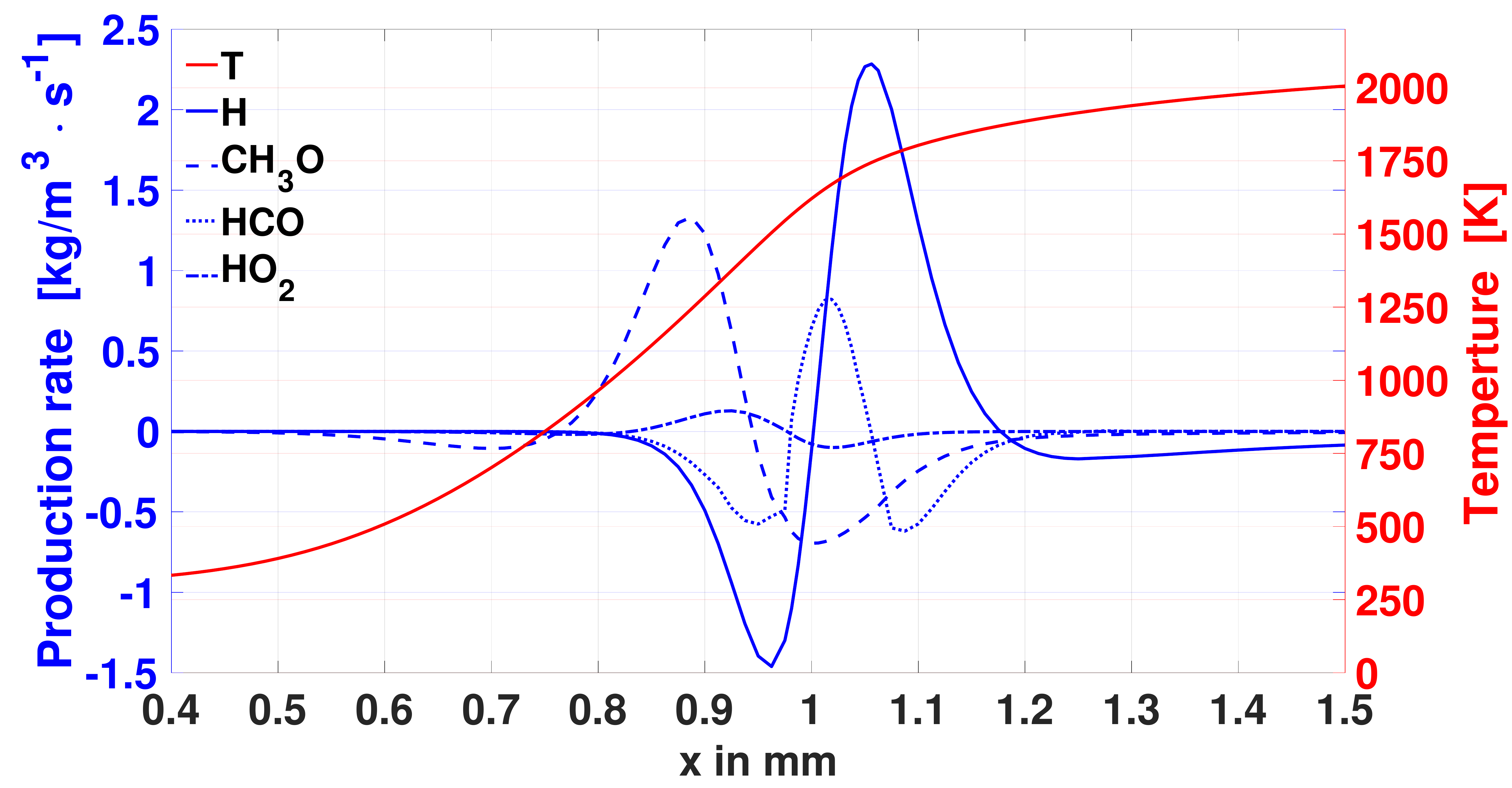}\\[2mm]
\parbox{13cm}{
\caption{One-dimensional study of the production rate of some
intermediate species. We choose HCO as candidate for the mesh design.}
\label{fig:methane-1D}
}
\end{figure}

\begin{figure}[ht!]
\centering
\includegraphics[width=0.9\linewidth]{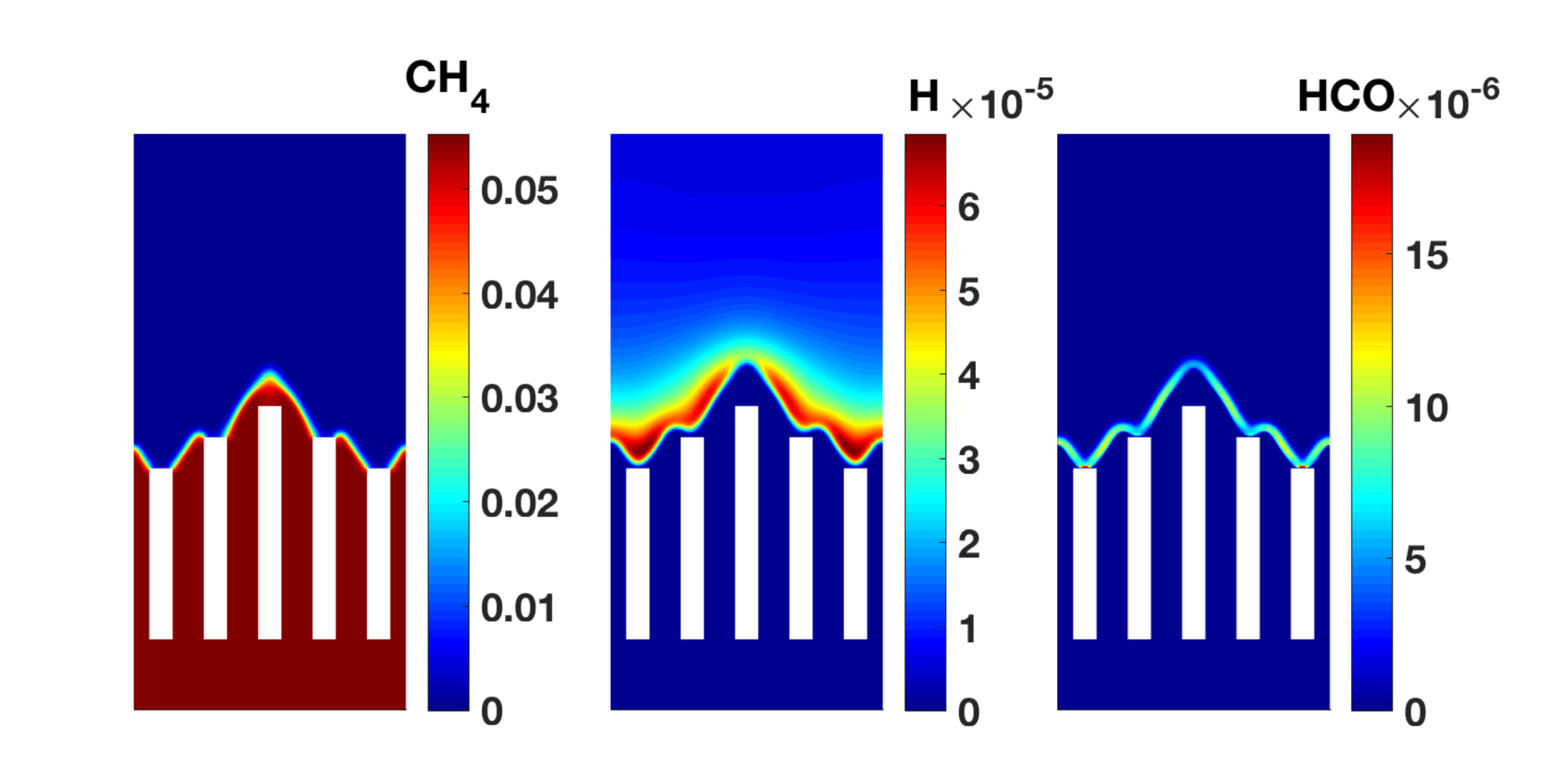}\\[2mm]
\parbox{13cm}{
\caption{Mass fraction of CH$_4$, H, and HCO at the steady state.}
\label{fig:methane-adapt-solution}
}
\end{figure}

\begin{figure}[ht]
\centering
\includegraphics[width=0.85\linewidth]{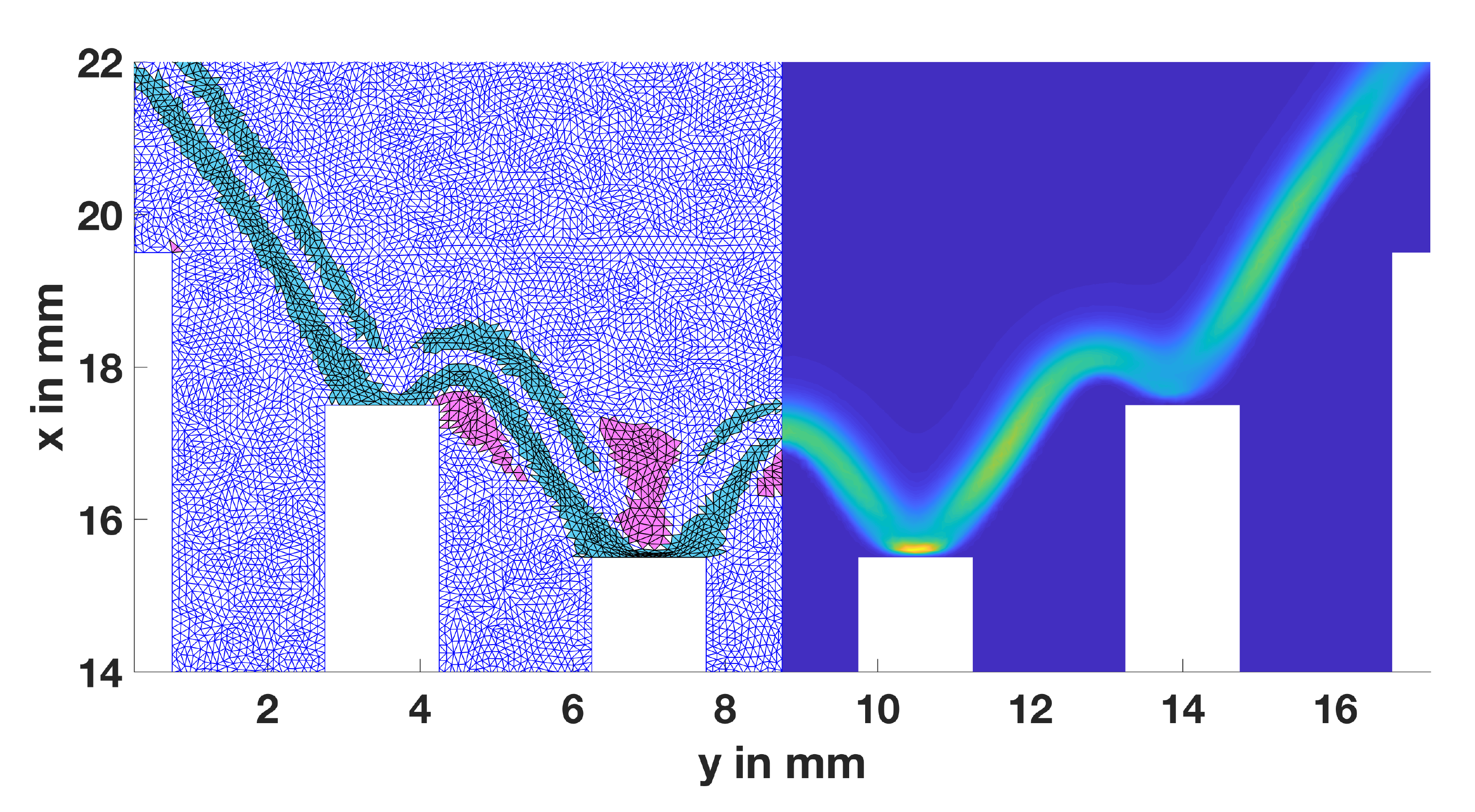}\\[2mm]
\parbox{13cm}{
\caption{Adaptive grid (left) and profile of the radical HCO (right) close to the top
of the lamellae at steady state. The elements colored in magenta show enlarged triangles
(up to $15\%$), and the elements colored in cyan illustrate the compressed
cells (up to $-5\%$). The mesh resolution is clearly aligned to the flame structure.}
\label{fig:methane-whco-adapt-mesh-sol}
}
\end{figure}

\begin{figure}[ht!]
\centering
\includegraphics[width=0.68\linewidth]{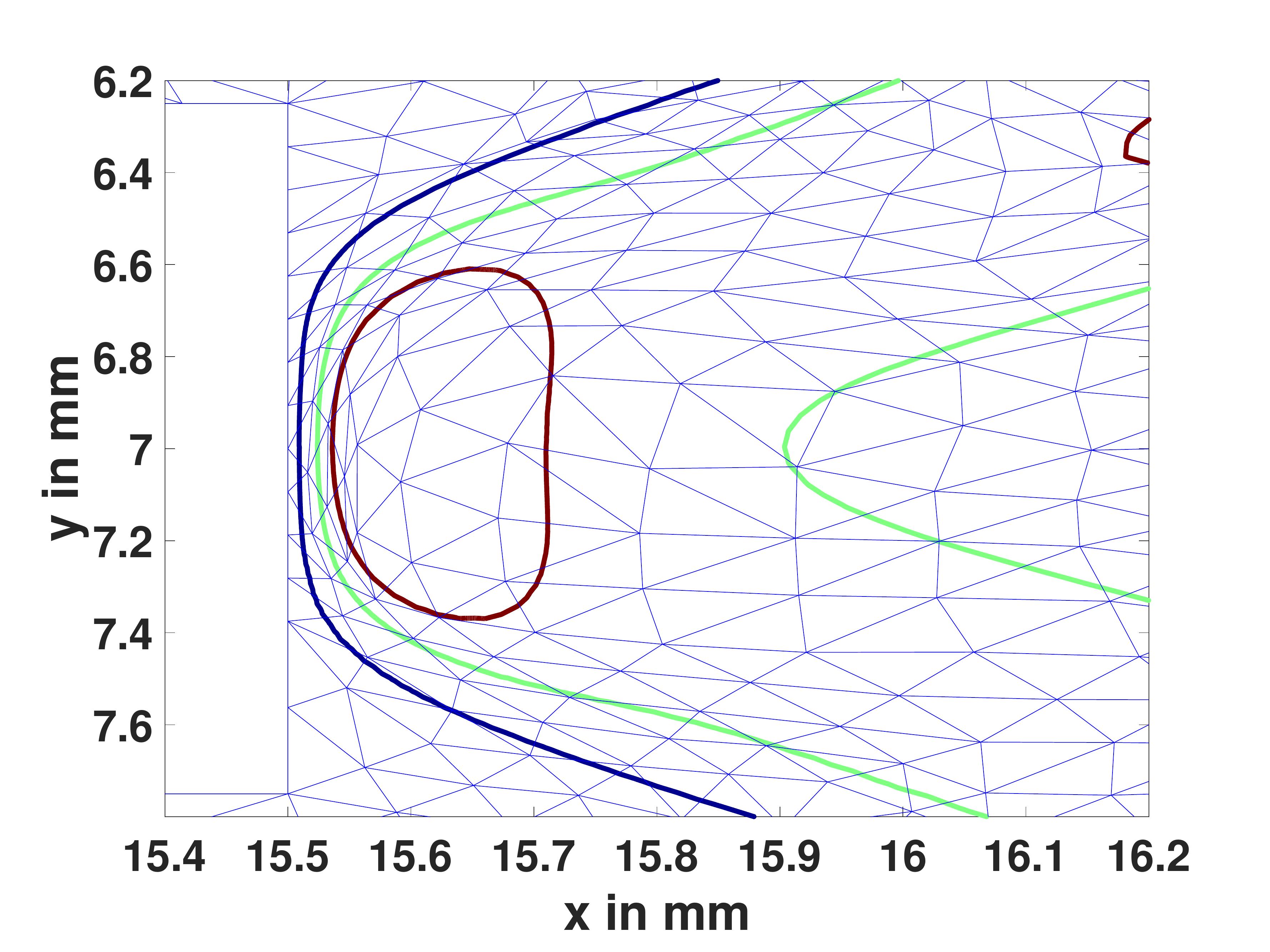}\\[2mm]
\parbox{13cm}{
\caption{Resolution of the flame layer of HCO at the top of the hottest
lamella. Colour curves represent
contour lines of the mass fraction of HCO: blue/green/brown correspond to
$0.1,\,0.3,\,0.5*\Vert w_{\text{HCO},h}\Vert_{\infty}$, respectively.
}
\label{fig:methane-zoom-adaptive-mesh}
}
\end{figure}

For the choice of the monitor function, we study a one-dimensional simplification
of the chemical reaction process. A corresponding simulation shows that the formyl
radical HCO and its production rate have an extremely thin flame layer, see
Fig.~\ref{fig:methane-1D}. Moreover, both are very sensible with respect to
insufficient mesh resolution \cite{Sun2018}. Hence, $w_{\text{HCO}}$ is a good candidate
for the mesh design. We set $\bpsi= D^2 w_{\text{HCO},h}$, which is now updated
in each time step. In Fig.~\ref{fig:methane-adapt-solution}, we show the mass fraction of
CH$_4$, H, and HCO at the steady state. The thin reaction zone of HCO is
clearly visible. A closer look at the top of the lamellae presented
in Fig.~\ref{fig:methane-whco-adapt-mesh-sol} and
Fig.~\ref{fig:methane-zoom-adaptive-mesh} reveals that
the mesh cells are strongly compressed there and aligned with respect
to the flame structure. Their areas are correspondingly reduced.

We would like to mention that the detection of the formyl radical
HCO is also of great interest since it can provide information
about the local heat release rate
which is a key parameter in the understanding of combustion processes.
Due to the low signal level, some optical diagnostic techniques, such
as HCO planar laser-induced fluorescence (PLIF), are not capable for
visualization and quantitative measurements
\cite{KieferLiSeegerLeipertzAlden2009}. Thus, the numerical investigation
of the dynamic behaviour of this radical can reveal
more information for further research.

\section{Conclusions}
In this paper, we have presented a two-dimensional adaptive moving mesh
method based on
combustion-specific design criteria in order to improve the numerical
resolution of steady state premixed flames in the low Mach number regime.
We have discussed the key ideas needed to design moving meshes which
allow a relocation of mesh points without changing their connectivity.
An iterative and simultaneous moving mesh strategy that balances solution
gradients or interpolation errors over the whole spatial mesh have been
applied to two combustion problems in 2D: a ozone decomposition and
a methane flame for a lamella burner. In both simulations, a
curvature-based adaptation of the mesh points yielded an improved flame
resolution and a stable calculation of the steady state. Adaptive moving
meshes led to an enhanced simulation capability that has made it possible
to simulate realistic flames without using more sophisticated local mesh
refinement methods.

\section{Acknowledgements}
Z. Sun and J. Lang are supported by the Deutsche Forschungsgemeinschaft within
the Graduate School of Excellence Energy Science and Engineering (DFG GSC1070).

\bibliographystyle{plain}
\bibliography{bibmmcomb}

\end{document}